\input amstex
%\magnification=\magstep1
\def\1{\text{\bf 1}}
\def\Aut{\text{\rm Aut}\,}
\def\BH{{\Cal B}({\Cal H})}

\def\CB{{\Cal C}{\Cal B}}
\def\Cb{{\Cal C}_b}
\def\Fb{{\Cal F}_b}
\def\UCb{{\Cal U}{\Cal C}_b}
\def\ULCb{{\Cal L}\,{\Cal U}{\Cal C}_b}
\def\URCb{{\Cal R}\,{\Cal U}{\Cal C}_b}
\def\id{\text{\rm id}}

\def\Ker{\text{\rm Ker}\,}
\def\Ran{\text{\rm Ran}\,}
\def\Lie{\text{\bf L}\,}

\documentstyle{amsppt}
\topmatter
\title Amenability, completely bounded projections, 
dynamical systems and smooth orbits\endtitle
\author Daniel Belti\c t\u a and Bebe Prunaru\endauthor
%\dedicatory  \enddedicatory
\address Institute of Mathematics ``Simion Stoilow'' of the Romanian
Academy,
P.O. Box 1-764, RO-014700 Bucharest, Romania\endaddress
\email Daniel.Beltita\@imar.ro, Bebe.Prunaru\@imar.ro\endemail
%\date May 26, 2004\enddate
\abstract
We describe a general method to construct completely bounded
idempotent mappings
on operator spaces, 
starting from amenable semigroups of completely bounded mappings. 
We then explore several applications of that method 
to injective operator spaces,
fixed points of completely contractive mappings, 
Toeplitz operators, dynamical systems and similarity orbits of group representations. 
\endabstract
\leftheadtext{Daniel Belti\c t\u a and Bebe Prunaru}
\rightheadtext{Amenability and completely bounded projections}
\endtopmatter

\document

\head 1. Introduction\endhead

If an injective von Neumann algebra is acted on by an amenable group 
then the corresponding fixed point algebra is in turn injective 
(Theorem~3.16 in Chapter~XV in \cite{Ta03}). 
This fact turns out to play a key role in 
several proofs in the theory of operator algebras. 
To give only one example in this connection, 
see the proof that an injective von Neumann algebra of type~III 
is semidiscrete 
(\S 3 in Chapter~XV in \cite{Ta03}). 

In the present paper we investigate what versions the aforementioned fact 
might have in the more general framework of operator spaces 
(see Corollary~3.2 below). 
Our initial motivation was that it might be useful to have a very general setting 
where completely bounded projections are associated with actions of semigroups. 
With the general result at hand 
(see Theorem~3.1 below) 
we soon realized that a lot of seemingly unrelated structures in operator theory 
can now be understood in a unifying manner. 
Thus such different things as dynamical systems, 
generalized Toeplitz operators 
or homogeneous spaces of Lie groups can be looked at 
from a unique point of view. 

We should point out that the technique of averaging over amenable groups 
has a long history in functional analysis and related areas. 
Its applications range from representation theory of finite and compact groups 
(see the so-called Weyl's unitary trick) to 
ergodic theory 
(see \cite{Lu92}) 
and cohomology of von Neumann algebras 
(see the papers \cite{SS98} and \cite{SS04}). 
From this point of view, what we are doing in the present paper 
is to investigate the relationship between that technique 
and the idea of completely bounded map. 

The structure of the paper is as follows: 
In Section~2 we introduce the notion of operator $S$-space,  
which is roughly speaking an operator space $X$ equipped with a semigroup $S$ 
of completely bounded maps. 
To each $S$-invariant subspace $Y\subseteq X^*$ 
and any left invariant mean on the corresponding 
space of coefficients ${\Cal C}_{X,Y}(S)$ 
we associate 
a completely bounded mapping $X\to Y^*$ 
and we study some basic properties of that construction. 
In Section~3 we prove our main result on 
existence of completely bounded projections on 
fixed point subspaces (Theorem~3.1) and 
then we explore some of the consequences 
of that theorem (see Theorems 3.4--3.5 and 
Corollaries 3.6 through 3.10). 

\specialhead Preliminaries\endspecialhead 

Our 
basic references for the theory of operator spaces and completely bounded maps 
are the monographs \cite{ER00}, \cite{Pa02} and~\cite{BL04}.  
We shall now recall several basic facts 
that will be needed in the sequel. 
If $X$ is a vector space and $p,q\ge1$ 
then $M_{p,q}(X)$ is the space of all $p$ by $q$ matrices 
with entries in $X$ and $M_p(X)=M_{p,p}(X)$. 
If $X$ and $Y$ are vector spaces, $\varphi\colon X\to Y$ 
is a linear mapping and $n\ge1$ then $\varphi_n\colon M_n(X)\to M_n(Y)$ 
is defined by 
$\varphi_n([x_{ij}])=[\varphi(x_{ij})]$ 
for every $[x_{ij}]\in M_n(X)$. 

Let ${\Cal H}$ be a complex Hilbert space 
and $\BH$ the $C^*$-algebra of all bounded linear operators on ${\Cal H}$. 
Then $M_n(\BH)$ has a unique $C^*$-algebra norm $\Vert\cdot\Vert_n$ 
induced by its identification with 
${\Cal B}({\Cal H}^{(n)})$  where 
${\Cal H}^{(n)}$ is the orthogonal sum of $n$ copies of ${\Cal H}$. 
An operator space is a complex vector space $X$ endowed with 
a complete norm $\Vert\cdot\Vert_n$ on every space $M_n(X)$ 
and with the property that there exists a linear mapping 
$\varphi\colon X\to\BH$ for some Hilbert space ${\Cal H}$ 
such that
$\varphi_n\colon(M_n(X),\Vert\cdot\Vert_n)\to
(M_n(\BH),\Vert\cdot\Vert_n)$ 
is isometric for all $n\ge1$. 

Any closed subspace of ${\Cal B}({\Cal H})$ inherits 
a canonical structure of operator space. 
In particular this holds true for $C^*$-algebras. 
More precisely, if $A$ is a $C^*$-algebra and $n\ge1$ 
then $M_n(A)$ has a unique $C^*$-algebra norm
that is induced by an arbitrary faithful representation 
of $A$ on a Hilbert space. 
If $A$ and $B$ are $C^*$-algebras and $\varphi\colon A\to B$ 
is linear then $\varphi$ is said to be completely positive 
if $\varphi_n$ is a positive map for all $n\ge1$. 

If $X$ and $Y$ are operator spaces and $\varphi\colon X\to Y$ 
is linear then $\varphi$ is said to be completely bounded 
if 
$\Vert\varphi\Vert_{\text{cb}}
\mathop{=}\limits^{\text{def}}
\sup\{\Vert\varphi_n\Vert\mid n\ge1\}<\infty$, 
and completely contractive if $\Vert\varphi\Vert_{\text{cb}}\le1$.
Moreover $\varphi$ is completely isometric if 
$\varphi_n$ is isometric for all $n\ge1$.
For $X$ and $Y$ operator spaces the space $\CB(X,Y)$ 
of all completely bounded maps between $X$ and $Y$ 
is a Banach space when endowed with the norm $\Vert\cdot\Vert_{\text{cb}}$. 
Moreover it has an operator space structure 
given by the isomorphisms 
$M_n(\CB(X,Y))\simeq\CB(X,M_n(Y))$. 
We shall always denote $\CB(X,X)=\CB(X)$. 
When $Y={\Bbb C}$ the space $X^*=\CB(X,{\Bbb C})$ 
is called the operator space dual of $X$. 
Now let us consider the operator space $X^{**}=(X^*)^*$. 
Then it can be shown that the canonical injection 
$J\colon X\hookrightarrow X^{**}$ is a complete isometry. 

If $X$, $Y$ and $Z$ are operator spaces and $\varphi\colon X\times Y\to Z$ 
is a bilinear map then for all $p,q\ge1$ one denotes 
$$\varphi_{p;q}\colon M_p(X)\times M_q(Y)\to M_{pq}(Z),\quad 
\varphi_{p;q}([u_{ij}],[v_{kl}])=[\varphi(u_{ij},v_{kl})]_{(i,k),(j,l)}.
$$ 
Then $\varphi$ is said to be completely bounded if 
$\Vert\varphi\Vert_{\text{cb}}=\sup\{\Vert\varphi_{p;q}\Vert\mid p,q\ge1\}<\infty$. 
As in the case of completely bounded linear maps, 
the space $\CB(X\times Y,Z)$ of all completely bounded 
bilinear maps $\varphi\colon X\times Y\to Z$ has an operator space structure; 
see \cite{ER00} for details. 

We shall also use the operator space projective tensor product 
$X\widehat{\otimes}Y$ of two operator spaces $X$ and $Y$ 
(see \cite{ER00} for the precise definition). 
All we need to know is that it is an operator space structure 
on a certain completion 
of the algebraic tensor product $X\otimes Y$ 
so that for any operator space $Z$ there is a canonical complete isometry 
$\CB(X\widehat{\otimes}Y,Z)\simeq\CB(X\times Y,Z)$. 
Moreover these spaces are completely isometric to 
$\CB(X,\CB(Y,Z))$. 

For a given Banach space $X$ there are two distinguished operator space structures on $X$, 
the maximal operator space 
$\max X=(X,\{\Vert\cdot\Vert_{\max,n}\}_{n\ge1})$ 
and the minimal operator space 
$\min X=(X,\{\Vert\cdot\Vert_{\min,n}\}_{n\ge1})$ 
such that for any other operator space structure $(X,\{\Vert\cdot\Vert_{n}\}_{n\ge1})$ 
one has 
$\Vert\cdot\Vert_{\min,n}\le\Vert\cdot\Vert_{n}\le\Vert\cdot\Vert_{\max,n}$ 
for all $n\ge1$.  
It can be shown that $(\max X)^*\simeq\min(X^*)$ 
and $(\min X)^*\simeq\max(X^*)$ as operator spaces 
(see (3.3.13) and (3.3.15) in \cite{ER00}). 

An operator algebra is an associative algebra $A$ endowed with an operator space structure
so that there exists a completely isometric homomorphism $\varphi\colon A\to\BH$ 
for some Hilbert space ${\Cal H}$. 
If moreover $A$ is an operator space dual and $\varphi$ can be chosen 
so that it is, additionally, weak$^*$-continuous, then $A$ is said to be 
a dual operator algebra. 
(See Chapter~2 in \cite{BL04} for details.)

We now recall a few basic definitions in differential geometry that will be needed in 
Corollary~3.10. 
A good reference for the differential geometry of Banach manifolds and homogeneous spaces 
is \cite{Up85}. 
Let $M$ be a Hausdorff topological space. 
A local chart of $M$ is a homeomorphism $\varphi\colon U\to V$, 
where $U$ is an open subset of $M$ and $V$ is an open subset of some real
Banach space. 
A smooth atlas of $M$ is any family of local charts $\{\varphi_j\colon U_j\to V_j\}_{j\in J}$ 
such that $\bigcup\limits_{j\in J}U_j=M$ 
and $\varphi_j\circ\varphi_k^{-1}\colon\varphi_k(U_j\cap U_k)\to\varphi_j(U_j\cap U_k)$ 
is a smooth mapping 
(between open subsets of Banach spaces) 
whenever $U_j\cap U_k\ne\emptyset$. 
A Banach manifold is a topological space $M$ equipped with a maximal smooth atlas. 
If $\widetilde{M}$ is another Banach manifold with a smooth atlas 
$\{\widetilde{\varphi}_{\widetilde{j}}\colon\widetilde{U}_{\widetilde{j}}\to
\widetilde{V}_{\widetilde{j}}\}_{\widetilde{j}\in\widetilde{J}}$ 
then a continuous mapping $f\colon M\to\widetilde{M}$ is smooth if 
$\widetilde{\varphi}_{\widetilde{j}}\circ f\circ\varphi_j^{-1}
\colon\varphi_j(f^{-1}(\widetilde{V}_{\widetilde{j}}))\to V_j$ 
is smooth 
whenever $f^{-1}(\widetilde{V}_{\widetilde{j}})\cap V_j\ne\emptyset$. 

A Banach-Lie group is a group $G$ which is also a Banach manifold 
such that the group operations 
(i.e., multiplication and inversion) are smooth. 
For instance, 
if $A$ is a unital associative Banach algebra then its group 
of invertible elements, denoted by $A^\times$, 
is a Banach-Lie group. 
Now let $G$ be a Banach-Lie group and $H$ a subgroup of $G$. 
We say that $H$ is a Banach-Lie subgroup if there exists a local chart $\varphi\colon U\to V$ 
of $G$ such that $\varphi(U\cap H)=V\cap{\Cal W}$, 
where $U$ is an open neighborhood of $\1\in G$, 
$V$ is an open subset of the Banach space ${\Cal Z}$, 
and ${\Cal W}$ is a split subspace of ${\Cal Z}$ 
(that is, there exists a bounded linear operator $E\colon{\Cal Z}\to{\Cal Z}$ 
such that $E^2=E$ and $\Ran E={\Cal W}$). 
If this is the case, 
then $G/H$ with the quotient topology has a structure of Banach manifold  such that 
the natural projection 
$\pi\colon G\to G/H$ 
is smooth and has smooth local cross-sections 
on a neighborhood of each point of $G/H$ 
(see e.g., Theorem~8.19 and Corollary~8.3 in \cite{Up85}). 
In this case we say that $G/H$ is a homogeneous space of $G$, 
and the natural transitive action 
$$G\times G/H\to G/H,\quad (g_1,g_2H)\mapsto g_1g_2H$$
is a smooth mapping. 
Now assume that $G=A^\times$ for some unital associative Banach algebra $A$ 
and that $H$ is an algebraic subgroup of $G$ (of degree $\le d$) in the sense that 
there exist an integer $d\ge1$ and a family ${\Cal F}$ of 
polynomial functions on $A\times A$ of degree $\le d$ 
such that 
$$H=\{g\in A^\times\mid(\forall f\in{\Cal F})\quad f(g,g^{-1})=0\}.$$
Denote 
$$\Lie(H)=\{a\in A\mid(\forall t\in{\Bbb R})\quad \exp(ta)\in H\}.$$
Then $H$ with the topology inherited from $A$ is a Banach-Lie group 
and $\Lie(H)$ (the Lie algebra of $H$) is a closed subspace of $A$ 
such that $[a,b]:=ab-ba\in\Lie(H)$ whenever $a,b\in\Lie(H)$. 
If it happens that $\Lie(H)$ is a split subspace of $A$, 
then $H$ is a Banach-Lie subgroup of $A^\times$. 
(See the main theorem in \cite{HK77} or Theorem~7.14 in \cite{Up85}.)

\head 2. Operator $S$-spaces\endhead

We begin this section by introducing some terminology on 
semitopological semigroups; we refer to \cite{BH67} and \cite{BJM78} 
for more details. 

\definition{Definition 2.1} 
For any semigroup $S$ we denote by $\Fb(S)$ 
the commutative unital $C^*$-algebra of all 
complex bounded functions on $S$ with the $\sup$ norm $\Vert\cdot\Vert_\infty$.
For each $t\in S$ we define 
$$L_t\colon\Fb(S)\to\Fb(S)\quad\text{and}\quad
R_t\colon\Fb(S)\to\Fb(S)$$ 
by 
$(L_tf)(s)=f(ts)$ and $(R_tf)(s)=f(st)$ 
whenever $s\in S$ and $f\in\Fb(S)$. 

Now assume that the semigroup $S$ is equipped with a topology. 
We say that $S$ is a {\it right} (respectively, {\it left}) 
{\it topological semigroup} 
if for each $s\in S$ the mapping $S\to S$, $t\mapsto ts$ 
(respectively, $t\mapsto st$) is continuous. 
Moreover $S$ is a {\it semitopological semigroup} if 
it is both left and right topological. 

If the semigroup $S$ is equipped with a topology then 
we denote by $\Cb(S)$ the set of all continuous functions in $\Fb(S)$.  
When $S$ is a right topological semigroup we denote 
$\ULCb(S)$ the set of all {\it left uniformly continuous} 
bounded complex functions on $S$. 
That is, 
$f\in\ULCb(S)$ if and only if 
$f\in\Cb(S)$ and the mapping 
$S\to\Cb(S)$, $s\mapsto R_sf$,
is continuous. 
Similarly, when $S$ is a left topological semigroup we define the set $\URCb(S)$ 
of all {\it right uniformly continuous} bounded complex functions on $S$ 
by the above condition with $R_s$ replaced by $L_s$. 
Moreover, when $S$ is a semitopological semigroup we shall need the set 
$\UCb(S)=\ULCb(S)\cap\URCb(S)$
consisting of all {\it uniformly continuous} bounded
complex functions on $S$. 
It is clear that all of the sets $\ULCb(S)$, $\URCb(S)$ and $\UCb(S)$ 
are unital $C^*$-subalgebras of $\Cb(S)$.  

Next assume again that $S$ is an arbitrary semigroup and 
let ${\Cal T}$ be any linear subspace of $\Fb(S)$. 
We say that ${\Cal T}$ is {\it unital} if it contains the unit element $\1$ 
of $\Fb(S)$ 
(i.e., if each constant function belongs to ${\Cal T}$). 
In this case, a {\it state} of ${\Cal T}$ 
is a linear functional $\mu\colon{\Cal T}\to{\Bbb C}$ 
such that $\Vert\mu\Vert=\mu(\1)=1$. 
Now assume that ${\Cal T}$ is a linear subspace of $\Fb(S)$ that is invariant under 
the operators $L_t$ for each $t\in S$.
We say that a linear functional $\mu\colon{\Cal T}\to{\Bbb C}$ is 
{\it $S$-invariant} if $\mu\circ L_t=\mu$ for all $t\in S$.  
The unital subspace ${\Cal T}$ of $\Fb(S)$ is said to be {\it amenable} 
if it admits an $S$-invariant state. 
If the space $\Fb(S)$ is 
amenable, then the semigroup $S$ is said to be {\it amenable}. 

A topological group  $S$ is said to be {\it amenable} 
if the space $\URCb(S)$ is amenable. 
For instance the unitary groups of all injective von Neumann algebras 
with the strong operator topology, 
and also the unitary groups of all nuclear unital $C^*$-algbras 
with the weak topology are amenable groups 
(see \cite{dlH79} and \cite{Pat92}). 
It is known that if $S$ is an amenable {\it locally compact} group 
then even the larger space 
$\Cb(S)$ is amenable
(see Theorem~2.2.1 in \cite{Gr69}). 
\qed
\enddefinition

\definition{Definition 2.2}
Let $S$ be a semigroup and $X$ an operator space. 
We say that $X$ is an {\it operator $S$-space}
if it
is equipped with a mapping 
$$\alpha\colon S\times X\to X,\quad (s,x)\mapsto
\alpha(s,x)=\alpha_s(x)$$
satisfying the following conditions: 

\itemitem{\rm(i)} for all $s,t\in S$ we have $\alpha_{st}=\alpha_s\circ\alpha_t$; 

\itemitem{\rm(ii)} for all $s\in S$ the mapping $\alpha_s\colon X\to X$ 
is completely bounded linear, and moreover 
$\sup\limits_{s\in 
S}\|\alpha_s\|_{\text{cb}}<\infty$.  

\noindent We say that $X$ is a {\it dual operator $S$-space} if moreover 
there exists an operator space $X_*$ such that $X=(X_*)^*$ and 
$(\alpha_s)^*X_*\subseteq X_*$ ($\subseteq X^*$) 
for all $s\in S$. 
An equivalent condition is that $\alpha_s\colon X\to X$ 
is weak$^*$-continuous for all $s\in S$.  
In this case, $X_*$ is said to be a {\it predual of the operator $S$-space} X. 

Let $X$ be an operator $S$-space and let $Y\subseteq X^*$ be a closed 
linear subspace such that $(\alpha_s)^*Y\subseteq Y$ for all $\in S$. 
We denote by ${\Cal C}_{X,Y}(S)$ 
the smallest unital closed subspace of $\Fb(S)$ 
that contains 
all the functions $f_{x,\psi}:=\psi(\alpha(\cdot,x))$ 
for $x\in X$ and $\psi\in Y$. 
We always think of ${\Cal C}_{X,Y}(S)$ as an operator space 
with the unique operator space structure that makes the inclusion map
${\Cal C}_{X,Y}(S)\hookrightarrow\Fb(S)$ into a complete isometry. 

Note that for all $s\in S$, $x\in X$ and $\psi\in Y$ 
we have $L_s(f_{x,\psi})=f_{x,(\alpha_s)^*\psi}$.  
Thus ${\Cal C}_{X,Y}(S)$ is invariant under $L_s$ for all $s\in S$. 
\qed
\enddefinition

\proclaim{Lemma 2.3}
Let $S$ be a semigroup, $X$ an operator $S$-space, 
$Y\subseteq X^*$ a closed linear subspace such that 
$(\alpha_s)^*Y\subseteq Y$ for all $\in S$, 
and ${\Cal C}_{X,Y}(S)$ as above. 
Then the mapping 
$$E^0\colon Y\times X\to{\Cal C}_{X,Y}(S),\quad 
(\psi,x)\mapsto(\psi\circ\alpha)(\cdot,x)$$
is a completely bounded bilinear mapping and 
$\|E^0\|_{\text{cb}}\le\sup\limits_{s\in S}\|\alpha_s\|_{\text{cb}}$. 
\endproclaim

\demo{Proof}
Let $p$, $q$ be arbitrary integers, denote ${\frak p}=\{1,2,\dots,p\}$ and 
${\frak q}=\{1,2,\dots,q\}$, 
and 
consider the bilinear mapping 
$$(E^0)_{p;q}\colon M_p(Y)\times M_q(X)\to M_{pq}(\Fb(S))$$
defined by 
$$(E^0)_{p;q}(\psi,x)=\bigl[\psi_{ij}(\alpha(\cdot,x_{kl}))
\bigr]_{(i,k),(j,l)\in{\frak p}\times{\frak q}}$$
for $\psi=(\psi_{i,j})_{i,j\in{\frak p}}\in M_p(Y)\subseteq M_p(X^*)\simeq\CB(X,M_p)$ 
and $x=(x_{kl})_{k,l\in{\frak q}}\in M_q(X)$. 
What we have to prove is that 
the norm of the bilinear mapping $(E^0)_{p;q}$ 
is at most 
$\sup\limits_{s\in S}\|\alpha_s\|_{\text{cb}}$. 
In fact, 
$$\aligned
\|(E^0)_{p;q}(\psi,x)\|
&=\sup_{s\in S}\|\bigl[\psi_{ij}(\alpha(\cdot,x_{kl}))
\bigr]_{(i,k),(j,l)\in{\frak p}\times{\frak q}}\| \cr
&=\sup_{s\in S}\|\psi_q((\alpha_s)_q(x))\| \quad\qquad\qquad\qquad
  (\text{compare (1.1.30) in \cite{ER00}})\cr
&\le\sup_{s\in S}\|\psi\|_{\text{cb}}\cdot\|\alpha_s\|\cdot\|x\| \cr
&=\sup\limits_{s\in S}\|\alpha_s\|_{\text{cb}} \cdot\|\psi\|\cdot\|x\|   
\quad\qquad\qquad\quad
  (\text{see (3.2.5) in \cite{ER00}})
\endaligned$$
and we are done.
\qed
\enddemo

\definition{Definition 2.4}
Let $S$ be a semigroup, $X$ an operator $S$-space with 
the semigroup action $\alpha\colon S\times X\to X$, 
and $Y\subseteq X^*$ a closed linear subspace such that 
$(\alpha_s)^*Y\subseteq Y$ for all $\in S$. 
Consider the bilinear map 
$$E^0\colon Y\times X\to{\Cal C}_{X,Y}(S),\quad 
(\psi,x)\mapsto(\psi\circ\alpha)(\cdot,x)$$ 
from Lemma~2.3. 
Then for each bounded linear functional  
$\mu\colon{\Cal C}_{X,Y}(S)\to{\Bbb C}$ 
we define the mapping 
$$E_\mu\colon X\to Y^*,\quad 
(E_\mu(x))(\psi)=\mu(E(\psi,x))=\mu((\psi\circ\alpha)(\cdot,x))$$
for all $x\in X$ and $\psi\in Y$.
\qed 
\enddefinition

\proclaim{Lemma 2.5}
With the notation of {\rm Lemma~2.3} and 
{\rm Definition~2.4} the bilinear mapping $E^0$ 
gives rise to a completely bounded linear mapping
$E\colon Y\widehat{\otimes} X\to{\Cal C}_{X,Y}(S)$ 
such that $E(\psi\otimes x)=f_{x,\psi}$. 
Its dual is 
a completely bounded linear mapping 
$$E^*\colon{\Cal C}_{X,Y}(S)^*\to \CB(X,Y^*)$$ 
with $\|E^*\|_{\text{cb}}\le\sup\limits_{s\in 
S}\|\alpha_s\|_{\text{cb}}$ and 
$E^*(\mu)=E_\mu$ for all $\mu\in{\Cal C}_{X,Y}(S)^*$. 
In particular, for each $\mu\in{\Cal C}_{X,Y}(S)^*$ 
we have 
$\|E_\mu\|_{\text{cb}}\le\sup\limits_{s\in 
S}\|\alpha_s\|_{\text{cb}}\cdot\|\mu\|$.
\endproclaim

\demo{Proof}
Denote $M:=\sup\limits_{s\in S}\|\alpha_s\|_{\text{cb}}$. 
It follows by the above Lemma~2.3 along with Proposition~7.1.2 
in~\cite{ER00} 
that the bilinear mapping 
$E\colon Y\times X\to{\Cal C}_{X,Y}(S)$ 
naturally corresponds to a completely bounded mapping 
$E\colon Y\widehat{\otimes} X\to{\Cal C}_{X,Y}(S)$ with 
$\|E\|_{\text{cb}}\le M$. 
Consequently, 
for the mapping dual to $E\colon Y\widehat{\otimes} X\to{\Cal C}_{X,Y}(S)$ 
we have  
$$E^*\colon{\Cal C}_{X,Y}(S)^*\to(Y\widehat{\otimes} X)^*=\CB(X,Y^*)$$
(the last equality follows by Corollary~7.1.5 in \cite{ER00}) 
and $\|E^*\|_{\text{cb}}=\|E\|_{\text{cb}}\le M$. 

Moreover, we have by Corollary~2.2.3 in \cite{ER00} 
that any continuous linear functional 
$\mu\colon{\Cal C}_{X,Y}(S)\to{\Bbb C}$ is completely bounded and 
$\|\mu\|_{\text{cb}}=\|\mu\|$. 
Consequently 
$\|E_\mu\|_{\text{cb}}=\|E^*(\mu)\|_{\text{cb}}
\le\|E^*\|_{\text{cb}}\cdot\|\mu\|\le M\cdot\|\mu\|$,
and the proof is finished.
\qed
\enddemo

\proclaim{Lemma 2.6}
Let $S$ be a semigroup, $X$ an operator $S$-space with 
the semigroup action $\alpha\colon S\times X\to X$, 
and $Y\subseteq X^*$ a closed linear subspace such that 
$(\alpha_s)^*Y\subseteq Y$ for all $\in S$. 
Let $E^*\colon{\Cal C}_{X,Y}(S)^*\to \CB(X,Y^*)$ as in 
{\rm Lemma~2.5} 
and endow the operator space $\CB(X,Y^*)$ 
with the semigroup action
$$\gamma\colon S\times\CB(X,Y^*)\to\CB(X,Y^*),\quad 
(s,\theta)\mapsto\gamma(t,\theta)=\gamma_t(\theta)
:=((\alpha_t)^*|_Y)^*\circ\theta.$$ 
Then for all $t\in S$ the diagram 
$$\CD
\CB(X,Y^*) @>{\gamma_t}>> \CB(X,Y^*) \cr
@A{E^*}AA @AA{E^*}A \cr
{\Cal C}_{X,Y}(S)^* @>{L_t^*}>> {\Cal C}_{X,Y}(S)^* 
\endCD$$
is commutative. 
In particular, 
if we have a bounded linear functional 
$\mu\colon{\Cal C}_{X,Y}(S)\to{\Bbb C}$ and an 
element $t\in S$ 
satisfying $\mu\circ L_t=\mu$, 
then $((\alpha_t)^*|_Y)^*\circ E_\mu=E_\mu$. 
\endproclaim

\demo{Proof}
Let $x\in X$ and $\psi\in Y$ arbitrary. 
We have 
$$\aligned
\langle((\alpha_t)^*|_Y)^*(E_\mu(x)),\psi\rangle
&=\langle E_\mu(x),(\alpha_t)^*(\psi)\rangle \cr
&=\langle E_\mu(x),\psi\circ\alpha_t\rangle \cr 
&=\mu\bigl((\psi\circ\alpha_t)(\alpha(\cdot,x))\bigr) \cr
&=\mu\bigl(\psi(\alpha(t\cdot,x))\bigr) 
   \qquad\qquad\qquad(\text{since $\alpha_t\alpha_s=\alpha_{ts}$})\cr
&=\mu\bigl(L_t((\psi\circ\alpha)(\cdot,x))\bigr) \cr
&=L_t^*\mu\bigl((\psi\circ\alpha)(\cdot,x)\bigr) \cr
&=\langle E_{L_t^*\mu}(x),\psi\rangle,
\endaligned$$
and the proof is complete.
\qed
\enddemo

\definition{Notation 2.7}
Let $S$ be a semigroup and $X$ an operator space such that there is a 
semigroup action 
$\alpha\colon S\times X\to X$, $(s,x)\mapsto\alpha(s,x)=\alpha_s(x)$, 
with $\alpha_s\in\CB(X)$ for all $s\in S$. 
Assume that $Y\subseteq X^*$ is a closed linear subspace such that 
$(\alpha_s)^*Y\subseteq Y$ for all $\in S$. 
Then there exists a natural semigroup action 
$$S\times Y^*\to Y^*,\quad 
(s,z)\mapsto((\alpha_s)^*|_Y)^*(z).$$
We denote 
$$X^S:=\{x\in X\mid(\forall s\in S)\quad \alpha_s(x)=x\},$$
and similarly 
$(Y^*)^S=\{z\in Y^*\mid
(\forall s\in S)\quad((\alpha_s)^*|_Y)^*(z)=z\}$. 

Also we denote by $\iota_Y\colon X\to Y^*$ the mapping 
defined by $\bigl(\iota_Y(x)\bigr)(y)=\langle y,x\rangle$ 
for all $x\in X$ and $y\in Y\subseteq X^*$. 
Thus $\iota_Y$ is the composition between the natural embedding 
$X\hookrightarrow X^{**}$ and 
the quotient mapping $X^{**}\to Y^*$, $\psi\mapsto \psi|_Y$. 
\qed
\enddefinition

\proclaim{Proposition 2.8}
Let $S$ be a semigroup, $X$ an operator $S$-space with the semigroup action 
$\alpha\colon S\times X\to X$, 
and $Y\subseteq X^*$ a closed linear subspace such that 
$(\alpha_s)^*Y\subseteq Y$ for all $\in S$. 
Assume that the space ${\Cal C}_{X,Y}(S)$ is amenable and pick 
an $S$-invariant state $\mu\in{\Cal C}_{X,Y}(S)^*$. 
Then 

\itemitem{\rm(i)} for all $x\in X^S$ we have $E_\mu(x)=\iota_Y(x)$, and 

\itemitem{\rm(ii)} $\iota_Y(X^S)\subseteq\Ran E_\mu\subseteq(Y^*)^S$.

\noindent In particular, if $X$ is a dual operator $S$-space 
and $Y^*=X$, then $\iota_Y=\id_X$, 
therefore $E_\mu(x)=x$ for all $x\in X^S$ and $\Ran E_\mu=X^S$. 
On the other hand, if $Y=X^*$ then $\iota_Y$ coincides with the canonical 
embedding $X\hookrightarrow X^{**}$ and, 
by this identification, 
it follows again that $E_\mu(x)=x$ for all $x\in X^S$ and 
$X^S\subseteq\Ran E_\mu\subseteq(X^{**})^S$. 
\endproclaim

\demo{Proof}
Assertion~(i) follows at once in view of the way $E_\mu$ was defined 
(see Definition~2.4) along with the fact that $\mu(\1)=1$, 
where $\1\in\Fb(S)$ is the function that is constant 1 on $S$. 

The first inclusion in assertion~(ii) follows by (i). 
The second inclusion follows by $S$-invariance of $\mu$ along with 
Lemma~2.6. 
\qed
\enddemo

In the following example we show that in the setting of Proposition~2.8 
it could happen that $X^S\ne\Ran E_\mu$. 

\example{Example 2.9}
Let $G$ be an amenable discrete infinite group, 
so that $\Cb(G)=\ell^\infty(G)$ 
is the $C^*$-algebra of all bounded complex functions on $G$, 
and there exists a $G$-invariant state $\mu\colon\ell^\infty(G)\to{\Bbb C}$. 
Consider the Banach space of all absolutely summable complex functions on $G$, 
$$\ell^1(G)=\Bigl\{f\colon G\to{\Bbb C}\mid 
\|f\|:=\sup_{\scriptstyle F\subseteq G\atop\scriptstyle F\text{ finite}}
\sum_{g\in F}|f(g)|<\infty\Bigr\},$$
and for all $g\in G$ define $\alpha_g\colon\ell^1(G)\to\ell^1(G)$ 
by $(\alpha_gf)(h)=f(g^{-1}h)$ for $h\in G$ and $f\in\ell^1(G)$. 
Also for each $f\in\ell^1(G)$ consider the convolution operator 
$$C_f\colon\ell^\infty(G)\to\ell^\infty(G), 
(C_fb)(h)=\sum_{g\in G}b(g)f(g^{-1}h)
\quad\text{for }b\in\ell^\infty(G)
\text{ and }h\in G.$$
Now define the operator $G$-space $X=\max\ell^1(G)$ 
(see Section~3.3 in \cite{ER00} for the definition of the 
functors $\min$ and $\max$ from Banach spaces to operator spaces) 
with 
$$\alpha\colon G\times X\to X,\quad (g,f)\mapsto\alpha_gf.$$
Then $X^*=\min\ell^\infty(G)$ and $X^{**}=\max(\ell^\infty(G))^*$ 
as operator spaces 
(see (3.3.13) and~(3.3.15) in~\cite{ER00}). 

Since $G$ is infinite, it follows that the only absolutely summable 
constant function on $G$ is $0$, hence 
$$X^G=\{0\}\ne(X^{**})^G,$$
where $(X^{**})^G$ is just the set of all 
$G$-invariant continuous linear functionals 
on the commutative $C^*$-algebra $\ell^\infty(G)$, 
and this set is different from $\{0\}$ since $G$ is amenable. 
Moreover, with the notation of Proposition~2.8 we claim that actually 
$$X^G=\{0\}\ne\Ran E_\mu\subseteq(X^{**})^G.$$
In fact it is easy to see that the mapping 
$E_\mu\colon\ell^1(G)\to(\ell^\infty(G))^*$ can be equivalently defined 
in terms of the convolution operators by 
$$(E_\mu f)(b)=\mu(C_fb) \quad\text{for }b\in\ell^\infty(G)
\text{ and }f\in\ell^1(G).$$
Hence for $f=\delta_{\1}$ (the characteristic function of $\{\1\}\subseteq G$) 
we have $E_\mu \delta_{\1}=\mu$, 
whence $0\ne\mu\in\Ran E_\mu$, 
and the above claim is proved. 
\endexample

\head 3. The main results\endhead

\proclaim{Theorem 3.1}
Let $S$ be a semigroup and $X$ an operator $S$-space with 
the semigroup action $\alpha\colon S\times X\to X$. 
Assume that one of the following hypotheses holds:

\itemitem{\rm(a)} $X$ is a dual operator $S$-space and 
$S$ is amenable as a discrete semigroup, or

\itemitem{\rm(b)} $X$ is a dual operator $S$-space, 
$S$ is an amenable locally compact topological group and 
for each $x\in X$ the mapping $\alpha(\cdot,x)$ 
is continuous with respect to the weak$^*$-topology of $X$, or 

\itemitem{\rm(c)} $X$ is a dual operator $S$-space, 
$S$ is an amenable topological group and the 
natural action of $S$ on the predual of $X$ is strongly continuous, 
or 

\itemitem{\rm(d)} $X$ is an operator $S$-space, 
$S$ is a compact left topological semigroup 
which is amenable as a discrete semigroup, and 
for each $x\in X$ the mapping $\alpha(\cdot,x)$ 
is continuous, or 

\itemitem{\rm(e)} $X$ is an operator $S$-space 
which is separable as a Banach space, 
$S$ is a compact left topological semigroup 
which is amenable as a discrete semigroup, and 
for each $x\in X$ the mapping $\alpha(\cdot,x)$ 
is weakly continuous, or 

\itemitem{\rm(f)} $S$ is a compact topological group 
and for each $x\in X$ the mapping $\alpha(\cdot,x)$ 
is weakly continuous. 
 
\noindent 
Then there exists a linear map $P\colon X\to X$ with the following properties: 

\itemitem{\rm(i)} $P\in\CB(X)$ and 
$\|P\|_{\text{cb}}\le\sup\limits_{s\in S}\|\alpha_s\|_{\text{cb}}$, 

\itemitem{\rm(ii)} $\Ran P=X^S$, and 

\itemitem{\rm(iii)} $P\circ P=P$. 

\endproclaim

\demo{Proof}
We first consider the conditions (a)--(c). 
Let $X_*$ be an operator space predual of $X$ as in Definition~2.2. 
We are going to make use of Proposition~2.8 for $Y=X_*$. 
To this end we first make sure that, if either of the conditions (a)--(c) 
is satisfied, then the function space ${\Cal C}_{X,Y}(S)$ is amenable. 
In the case (a), this is obvious since the space $\Fb(S)$ of all 
bounded complex functions on $S$ is amenable. 
In the case (b), we have ${\Cal C}_{X,Y}(S)\subseteq \Cb(S)$ 
and the space $\Cb(S)$ is amenable since the group $S$ is locally compact.
Finally, in the case (c) recall that 
for all $s\in S$, $x\in X$ and $\psi\in Y$ 
we have $L_s(f_{x,\psi})=f_{x,(\alpha_s)^*\psi}$. 
The hypothesis~(c) means that the mapping $S\to Y$, $s\mapsto(\alpha_s)^*\psi$,  
is continuous for all $\psi\in Y$, hence we get 
${\Cal C}_{X,Y}(S)\subseteq\URCb(S)$, 
while the latter space is amenable since $G$ is an amenable group. 
Consequently the space ${\Cal C}_{X,Y}(S)$ is amenable in either of 
the cases (a)--(c). 
Now pick an $S$-invariant state $\mu\colon{\Cal C}_{X,Y}(S)\to{\Bbb C}$ 
and denote $P=E_\mu\colon X\to Y^*=X$. 
Then we have by Proposition~2.8 that $\Ran P=X^S$ and $P$ is the identity map on $X^S$, 
whence the desired properties (ii)--(iii) follow. 
As for property (i), it is a consequence of Lemma~2.5. 

We now address the conditions (d)--(f). 
We are going to apply Proposition~2.8 with $Y=X^*$. 
Again we first need to check that the space ${\Cal C}_{X,Y}(S)$ is amenable. 
In both cases (d) and (e) 
this is obvious since the whole space $\Fb(S)$ is amenable. 
In the case (f) note that ${\Cal C}_{X,Y}\subseteq\Cb(S)$ 
and the latter space is amenable. 
Thus the space ${\Cal C}_{X,Y}(S)$ is amenable under either of the conditions (d)--(f), 
and then we can pick an $S$-invariant state $\mu\colon{\Cal C}_{X,Y}(S)\to{\Bbb C}$ 
and denote $P=E_\mu\colon X\to Y^*=X^{**}$. 
We are going to prove that actually $\Ran E_\mu\subseteq X$, 
and then the desired properties (i)--(iii) will follow just as above, 
by Proposition~2.8 along with Lemma~2.5. 

Firstly assume that the condition (d) is satisfied and 
let $x\in X$ arbitrary. 
In order to show that $E_\mu(x)\in X$ 
we have to check that $E_\mu(x)\colon X^*\to{\Bbb C}$ 
is weak$^*$-continuous. 
To this end, it is enough to check that $E_\mu(x)$ is weak$^*$-continuous 
on the unit ball of $X^*$. 
(See e.g., Corollary~2 to Theorem~6.2 in Chapter~IV 
of \cite{Sch66}.) 
Thus let $\{\psi_j\}_{j\in J}$ be a net in $X^*$ 
such that $\|\psi_j\|\le1$ for all $j\in J$ 
and $\psi_j\mathop{\longrightarrow}\limits^{\text{weak$^*$}}_{j\in J}0$. 
Then $\psi_j\mathop{\longrightarrow}\limits_{j\in J}0$ 
uniformly on the compact subsets of $X$. 
On the other hand, since $S$ is compact,  
it follows that $\{\alpha(s,x)\mid s\in S\}$ is a compact subset 
of $X$, hence 
$(\psi_j\circ\alpha)(\cdot,x)\mathop{\longrightarrow}\limits_{j\in J}0$ 
uniformly on $S$. 
Consequently 
$(E_\mu(x))(\psi_j)=\mu((\psi_j\circ\alpha)(\cdot,x))
\mathop{\longrightarrow}\limits^{\text{weak$^*$}}_{j\in J}0$, 
and thus $E_\mu(x)$ is weak$^*$-continuous
on the unit ball of $X^*$.

In the case (e), first recall that $\mu$ actually extends to 
an $S$-invariant state of $\Fb(S)$, and in particular 
to an $S$-invariant state $\mu\colon\Cb(S)\to{\Bbb C}$. 
Thus $\mu$ actually defines a Radon measure on $S$. 
Next, as in the case (d), we let $x\in X$ arbitrary and check that 
$E_\mu(x)\colon X\to{\Bbb C}$ is weak$^*$-continuous on the unit ball of $X^*$. 
Since $X$ is separable, the weak$^*$-topology of 
the unit ball of $X^*$ is metrizable, 
hence it is enough to check that, 
if $\{\psi_j\}_{j\ge0}$ is a {\it sequence} in the unit ball 
of $X^*$ with 
$\psi_j\mathop{\longrightarrow}\limits^{\text{weak$^*$}}_{j\to\infty}0$ 
then $\lim\limits_{j\to\infty}\bigl(E_\mu(x)\bigr)(\psi_j)=0$. 
But this fact follows by Lebesgue's dominated convergence theorem,  
since $\bigl(E_\mu(x)\bigr)(\psi_j)=\mu((\psi_j\circ\alpha)(\cdot,x))$ 
and $\Vert(\psi_j\circ\alpha)(\cdot,x)\Vert_\infty
\le\sup\limits_{s\in S}\Vert\alpha_s\Vert_{\text{cb}}\cdot\Vert x\Vert$ 
for all $j\ge1$. 

In the case (f), since $S$ is a compact group, 
it follows by Proposition~4.2.2.1 in \cite{Wa72} 
that 
for each $x\in X$ the mapping $\alpha(\cdot,x)$ is actually continuous, 
hence the conclusion follows by (d).  
Alternatively, note that 
the invariant state $\mu\colon\Cb(S)\to{\Bbb C}$ 
is defined by a probability Haar measure on $S$, 
and use of Proposition~2 and 
Remark~1 in Chapter~III, \S 4, no.~1 
in~\cite{Bo65} to show that $\Ran E_\mu(x)\in X$ for all $x\in X$. 
\qed
\enddemo

We note that, under the hypothesis (b) of Theorem~3.1, the mapping 
$E_\mu\colon X\to X$ used in the proof 
shows up in several places in the existing literature. 
See e.g., Weyl's unitary trick 
(that is, the fact that every representation of a compact group 
is similar to a unitary representation) 
or, more recently, 
\cite{DLRZ02}~and~\cite{OR03}. 

For the first corollary of Theorem~3.1 we recall that 
an operator space $Y$ is said to be {\it injective} 
if for any complete isometry $\varphi\colon X_0\to X$ 
and every $\psi_0\in\CB(X_0,Y)$ 
there exists $\psi\in\CB(X,Y)$ such that $\psi\circ\varphi=\psi_0$ 
and $\Vert\psi\Vert_{\text{cb}}=\Vert\psi_0\Vert_{\text{cb}}$. 
(See \cite{ER00} for details.) 

\proclaim{Corollary 3.2} 
Let $X$ be a dual operator space and $S$ an amenable semigroup 
of completely contractive, weak$^*$-continuous linear mappings on $X$. 
If $X$ is an injective operator space, then $X^S$ is in turn injective. 
\endproclaim

\demo{Proof}
It follows by Theorem~3.1 along with condition~(a) in Definition~2.2 
that there exists a completely contractive projection $P\colon X\to X$ 
with $\Ran P=X^S$. 
Now the desired conclusion follows by Proposition~4.1.6 in \cite{ER00}.  
\qed
\enddemo

It is safe to say that  most of the assertions 
contained in the following two theorems are parts of 
the folklore of operator algebras. 
However we would 
like to show how they follow directly from Theorem~3.1 
and to emphasize that the idempotent mappings we construct here 
are completely bounded. 
Before going further, we recall that any $*$-homomorphism of $C^*$-algebras 
is completely contractive and any $*$-automorphism 
of a von Neumann algebra is weak$^*$-continuous. 

\proclaim{Theorem 3.3}
Let $A$ be a $C^*$-algebra, $G$ a topological 
group, and $\alpha\colon G\to \Aut(A)$ a group 
homomorphism such that for each $x\in A$ the map 
$g\mapsto\alpha\sb g(x)$ is continuous with respect to the 
norm topology of $A$. 
(When $G$ is locally compact, 
the triple $(A,G,\alpha)$ with these properties is 
called in literature a $C^*$-dynamical system.) 
Then  the following hold true:

\itemitem{\rm(a)} 
If $G$ is amenable then there exists a completely 
bounded idempotent mapping $Q\colon A^*\to A^*$ whose range 
consists of all linear forms $\phi\in A^*$ which are 
$\alpha$-invariant, i.e., 
$\phi(\alpha_g(x))=\phi(x)$ for all $g\in G$ and all $x\in A$. 
Moreover, if $A$ is unital, then $Q(S(A))\subset 
S(A)$ hence $Q$ maps the set of all 
states of $A$ onto the set of all $\alpha$-invariant 
states. 

\itemitem{\rm(b)} 
If $G$ is a compact group, then there exists a 
completely positive and completely contractive 
idempotent $P\colon A\to A$ with 
$\Ran P=\{x\in A\mid\alpha_g(x)=x\text{ for all }g\in G\}$. 

\endproclaim

\demo{Proof}
(a) Let us consider the action $\beta$ of $G$ on the 
dual space $A^*$ defined by
$$\beta(g,\phi)(x)=\phi(\alpha\sb{g^{-1}}(x))$$
for all $x\in A$ and all $\phi\in A^*$. 
It easy to see 
that $A^*$ becomes, via this action, a dual operator 
$G$-space satisfying Theorem 3.1 item~(b) for the 
case when $G$ is locally compact or item~(c) for the 
case when $G$ is a topological group. Now the 
conclusion follows from that theorem.

b) Follows immediately from Theorem~3.1(f)
with $X=A$.
\qed
\enddemo 

\remark{Remark 3.4} 
We refer to \cite{Pe79} for general information of $C^*$-dynamical systems. 
In the case when $G$ is an amenable locally 
compact group, the item~(a) in Theorem~3.3 
holds true under the 
weaker hypothesis that all the functions 
$g\mapsto\phi(\alpha_g(x))$ are continuous on $G$ for all 
$x\in A$ and all $\phi\in A^*$. 
Indeed, in this case 
we can apply Theorem 3.1 item~(b) to the operator $G$-space $A$.
\qed 
\endremark

\proclaim{Theorem 3.5}
Let $M$ be a von Neumann algebra, $G$ a 
topological group and $\alpha\colon G\to \Aut(M)$ a group 
homomorphism such that for each $x\in M$ and each 
$\phi\in M_{*}$ (the predual of $M$) the functions 
$g\mapsto\phi(\alpha_g(x))$ are continuous on 
$G$. 
(When $G$ is locally compact, the triple 
$(M,G,\alpha)$ as above is called a $W^*$-dynamical system.) 
Then the following hold true:

\itemitem{\rm(a)} 
Suppose either  $G$ is an amenable locally 
compact group or $G$ is an amenable (not 
necessarily locally compact) topological group with 
the additional hypothesis (in this general case) that 
for all $\phi\in M_{*}$ the mapping 
$g\mapsto\alpha_g^*(\phi)$ is continuous 
with respect to the norm topology of the predual 
$M_{*}$. 
Then there exists a completely positive 
unital (hence completely contractive too) idempotent 
mapping $P\colon M\to M$ whose range is the fixed point 
algebra of $\alpha$. 
In particular, it follows that if 
$M$ is injective, then the fixed point algebra is also 
injective.    

\itemitem{\rm(b)} If $G$ is compact, then there exists a completely 
contractive idempotent mapping $Q\colon M_{*}\to M_{*}$ whose 
range is precisely the set of all $\alpha$-invariant 
normal forms on $M$. 
Moreover, $Q$ maps the set 
of all normal states of $M$ onto the set of all normal 
and $\alpha$-invariant states of $M$. 
The dual map
$Q^*\colon M\to M$ is a faithful completely positive  and 
normal idempotent mapping whose range is the fixed point 
algebra of $\alpha$ (faithful means that 
$\Ker P\cap M^{+}=\{0\}$).   

\endproclaim

\demo{Proof}
a) This follows from Theorem~3.1 item~(b) for the case 
when $G$ is locally compact or from item~(c) when 
$G$ is a topological group.

b) Let us consider, in a similar way as in the proof of 
the preceding theorem,  the action 
$$\beta\colon G\times M_{*}\to M_{*}$$
defined by
$$\beta(g,\phi)(x)=\phi(\alpha\sb{g^{-1}}(x)).$$
Then $M_{*}$ becomes an operator $G$-space, and 
moreover, by Proposition~4.2.2.1 in \cite{Wa72}, 
the action $\beta$ is also 
continuous with respect to the norm topology on $M_{*}$. 
Now the 
existence and other properties (except faithfulness) 
of $Q^*$ follows from item~(f) in Theorem~3.1. 
The expression of $Q$ as an integral with respect to the 
Haar measure on a compact group shows that $Q^*=P$, 
where $P$ is the one from~(a). 
Thence the asserted properties of $Q^*$ follow. 
\qed
\enddemo

As another consequence of Theorem~3.1 we now get the following 
version of Theorem~16(b) in \cite{Ke02}.  
See \cite{Ke04} for more information on generalized Toeplitz operators.

\proclaim{Corollary 3.6}
Let $(S,\cdot)$ be an amenable semigroup,  
${\Cal H}$ be a complex Hilbert space and $\rho\colon S\to{\Cal B}({\Cal H})$ 
a norm-continuous mapping such that 
$\rho(st)=\rho(s)\rho(t)$, $\rho(\1)=\id_{\Cal H}$ 
and $\|\rho(s)\|\le1$ for all $s,t\in S$. 
Now consider the space of $\rho$-Toeplitz operators 
$${\Cal T}(\rho)=\{C\in{\Cal B}({\Cal H})\mid(\forall s\in S)\quad 
\rho(s)C\rho(s)^*=C\}.$$
Then there exists a completely positive, completely contractive mapping 
$$P\colon{\Cal B}({\Cal H})\to{\Cal B}({\Cal H})$$ 
with 
$\Ran P={\Cal T}(\rho)$, $P\circ P=P$ and 
$P(ADB^*)=AP(D)B^*$ 
whenever $D\in{\Cal B}({\Cal H})$ and $A$ and $B$ belong to the commutant of 
$\rho(S)$.
\endproclaim

\demo{Proof}
First note that condition~(a) in Definition~2.2 is satisfied 
for $X={\Cal B}({\Cal H})$ 
with the structure of dual operator $S$-space defined by 
$$\alpha\colon S\times{\Cal B}({\Cal H})\to{\Cal B}({\Cal H}),\quad 
\alpha(s,A)=\rho(s)A\rho(s)^*.$$
Clearly ${\Cal B}({\Cal H})^S={\Cal T}(\rho)$, 
hence Theorem~3.1 shows that there exists an idempotent completely contractive 
linear mapping 
$P\colon{\Cal B}({\Cal H})\to{\Cal B}({\Cal H})$ with 
$\Ran P={\Cal T}(\rho)$. 
Now it follows by the very construction of $P$ 
 that $P$ is completely positive and 
$P(ADB^*)=AP(D)B^*$ 
whenever $D\in{\Cal B}({\Cal H})$ and $A$ and $B$ belong to the commutant of 
$\rho(S)$. 
\qed
\enddemo

\proclaim{Corollary 3.7}
Let ${\Cal T}({\Bbb T})$ be the space of all 
Toeplitz operators on the Hardy space $H^2({\Bbb T})$ associated with the unit disk. 
Then there exists a completely positive, completely contractive linear mapping 
$P\colon{\Cal B}(H^2({\Bbb T}))\to{\Cal B}(H^2({\Bbb T}))$ 
such that $P\circ P=P$ and $\Ran P={\Cal T}({\Bbb T})$. 
In particular ${\Cal T}({\Bbb T})$ is an injective operator space. 
\endproclaim

\demo{Proof}
Let 
$$M_z\colon H^2({\Bbb T})\to H^2({\Bbb T}),\quad 
(M_zf)(\text{e}^{\text{i}\theta})=\text{e}^{\text{i}\theta} f(\text{e}^{\text{i}\theta}),$$
the unilateral shift operator. 
It is well known that 
$${\Cal T}({\Bbb T})=\{C\in{\Cal B}(H^2({\Bbb T}))\mid M_z^*CM_z=C\},$$
hence the desired conclusion follows by Corollary~3.4 applied 
for the Abelian semigroup $(S,\cdot)=({\Bbb N},+)$ 
and $\rho\colon{\Bbb N}\to{\Cal B}(H^2({\Bbb T}))$, $\rho(n)C=(M_z^*)^n$ 
for all $n\in{\Bbb N}$. 
We note that $({\Bbb N},+)$ is amenable since it is Abelian
(see Theorem~1.2.1 in~\cite{Gr69}). 
\qed
\enddemo

As another consequence of Theorem~3.1 we now provide an alternative proof 
of Theorem~2.4(a) in \cite{AGG02}. 
In the special case when ${\Cal M}={\Cal B}({\Cal H})$, 
the next corollary shows that 
the set $C_{=}(\varphi)=\{X\in{\Cal B}({\Cal H})\mid\varphi(X)=X\}$ studied in 
\cite{Po03} 
is an injective operator space 
provided $\varphi\colon{\Cal B}({\Cal H})\to{\Cal B}({\Cal H})$ 
is a weak$^*$-continuous, completely positive, completely contractive map. 

\proclaim{Corollary 3.8}
Let ${\Cal M}$ be a $W^*$-algebra, $\alpha\colon{\Cal M}\to{\Cal M}$ 
a weak$^*$-continuous completely positive, completely contractive linear mapping, 
and denote 
$${\Cal M}^\alpha=\{x\in{\Cal M}\mid\alpha(x)=x\}.$$ 
Then there exists an idempotent, completely positive, completely contractive, linear mapping 
$P\colon{\Cal M}\to{\Cal M}$ 
with $\Ran P={\Cal M}^\alpha$. 
\endproclaim

\demo{Proof}
The existence of a completely contractive projection $P$ 
from ${\Cal M}$ 
onto ${\Cal M}^\alpha$  follows by 
Theorem~3.1 for $X={\Cal M}$, $S=({\Bbb N},+)$ 
and 
$$\alpha(n,x)=\alpha^n(x)$$ 
whenever $n\in{\Bbb N}$ and $x\in{\Cal M}$. 
To conclude the proof, we only have to remark that 
the idempotent mapping $P$ given by Theorem~3.1 
is completely positive according to its construction, since 
$\alpha\colon{\Cal M}\to{\Cal M}$ is completely positive 
(see also the construction of $E_\mu$ in Definition~2.4). 
\qed
\enddemo

We now arrive at a corollary that has interesting consequences in 
providing certain homogeneous spaces with structures of Banach manifolds. 
See \cite{BR04} and also Corollary~3.10 below. 

\proclaim{Corollary 3.9}
Let ${\Cal X}$ be a complex Banach space, $S$ a topological group 
and 
$$\alpha\colon S\to{\Cal B}({\Cal X}),\quad s\mapsto\alpha_s,$$
a norm continuous representation of $S$ by bounded linear operators on 
${\Cal X}$ 
such that $\alpha_{\1}=\id_{\Cal X}$ and 
$\sup\limits_{s\in S}\|\alpha_s\|<\infty$. 
Assume that one of the following hypotheses holds:

\itemitem{\rm(a)} $S$ is an amenable topological group and 
%there exists a complex Banach space ${\Cal C}$ 
%such that ${\Cal X}={\Cal C}^*$ and $(\alpha_s)^*{\Cal C}\subseteq{\Cal C}$ 
${\Cal X}$ is a dual Banach space such that $\alpha_s\colon{\Cal X}\to{\Cal X}$ 
is weak$^*$-continuous for all $s\in S$, or 

\itemitem{\rm(b)} 
%$S$ is a compact left topological semigroup 
%and ${\Cal C}_X(S)$ is amenable (for instance, 
$S$ is a compact topological group%)
.
%(for instance $S$ is a compact group and 
%the $S$-invariant state is defined by the probability Haar measure on $S$). 
 
\noindent Next denote
$${\Cal X}^S=\{x\in{\Cal X}\mid 
(\forall s\in S)\quad\alpha_s(x)=x\}.$$
Then there exists a bounded linear operator $P\in{\Cal B}({\Cal X})$ 
such that $\|P\|\le\sup\limits_{s\in S}\|\alpha_s\|$, 
$\Ran P={\Cal X}^S$ and $P^2=P$. 
\endproclaim

\demo{Proof}
We are going to apply Theorem~3.1 for the operator space $X=\max{\Cal X}$. 
According to~(3.3.9) in~\cite{ER00} we have an isometric identification 
${\Cal B}({\Cal X})\simeq\CB(\max{\Cal X})$, 
hence it follows at once that $\max{\Cal X}$ is an operator $S$-space. 
On the other hand, the above identification shows that the desired 
conclusion will 
follow as soon as we show that each of the present hypotheses implies 
one of the conditions of Theorem~3.1 for the operator space 
$X=\max{\Cal X}$. 

Actually, it is obvious that the present hypothesis~(b) implies 
that condition~(f) in Theorem~3.1 is fulfilled. 
As for the present hypothesis~(a), note that it implies that 
the condition~(c) in Theorem~3.1 is satisfied. 
In fact, it follows by (3.3.15) in \cite{ER00} that 
if ${\Cal Y}$ is a Banach space such that ${\Cal X}={\Cal Y}^*$ then 
$(\min{\Cal Y})^*=\max{\Cal Y}^*=\max{\Cal X}$, 
hence $\max{\Cal X}$ is the dual operator space of $\min{\Cal Y}$, and we 
are done. 
\qed
\enddemo

The next result is a partial extension of Theorem~4.8 in \cite{CG99} 
and is also related to Theorems 3.12~and~4.4 in \cite{ACS95}. 
We note that under hypothesis~(a) of this corollary we do {\it not} 
require that the group $G$ should be locally compact, 
and thus the result holds for infinite-dimensional Lie groups. 

\proclaim{Corollary 3.10}
Let $A$ be a unital operator algebra and denote by $A^\times$ 
its group of invertible elements. 
Consider an amenable topological group $G$ and denote 
$${\Cal R}:=\{\rho\colon G\to A^\times\mid 
\rho\text{ continuous group homomorphism and }\sup_{g\in 
G}\|\rho(g)\|<\infty\}.$$
Assume that one of the following conditions is satisfied: 

\itemitem{\rm(i)} $A$ is a dual algebra, or 

\itemitem{\rm(ii)} $G$ is compact.

\noindent Then the orbits of the action 
$$A^{\times}\times{\Cal R}\to{\Cal R},\quad (a,\rho)\mapsto a\cdot\rho(\cdot)\cdot 
a^{-1},$$
have natural structures of Banach manifolds that are smoothly acted 
on 
by the Banach-Lie group $A^\times$. 
\endproclaim

\demo{Proof}
In the proof we need  techniques and ideas from Lie theory 
that were recalled in the Introduction. 
Fix $\rho\in{\Cal R}$ and consider its isotropy group 
$$(A^\times)_\rho=\{a\in A^\times\mid(\forall g\in G)\quad 
a\cdot\rho(g)\cdot a^{-1}=\rho(g)\}.$$
We shall prove that $(A^\times)_\rho$ is a Banach-Lie subgroup of 
$A^\times$, and then 
the desired conclusion follows by Theorem~8.19 in \cite{Up85} 
in view of the natural bijection that exists from 
$A^\times/(A^\times)_\rho$ onto the orbit of $\rho$. 

To show that $(A^\times)_\rho$ is a Banach-Lie subgroup of $A^\times$, 
we first note that it is an algebraic subgroup of $A^\times$ of degree $\le1$ 
in the sense explained in the Introduction to the present paper. 
It then follows that $(A^\times)_\rho$ has a structure of Banach-Lie group with the topology 
inherited from 
$A^\times$, as a consequence of the main result of \cite{HK77}. 
The Lie algebra of $A^\times$ is 
$\Lie(A^\times)=A$,
while the Lie algebra of $(A^\times)_\rho$ is 
$$\Lie((A^\times)_\rho)=\{a\in A\mid(\forall g\in G)\quad 
a\cdot\rho(g)=\rho(g)\cdot a\}=\rho(G)',$$
hence it remains to prove that $\rho(G)'$ has a complement in $A$. 

To this end, consider the action of $G$ on $A$ defined by 
$$\alpha\colon G\times A\to A,\quad\alpha(g,a)=\rho(g)a\rho(g)^{-1}.$$
This action makes $A$ into an operator $G$-space, since $\alpha(g,\cdot)$ 
is completely bounded on $A$ for all $g\in G$ as an easy consequence of 
Theorem~17.1.2 in~\cite{ER00}. 
Moreover note that condition~(b) in Definition~2.2 is satisfied.
In case~(i), it follows by Theorem~2.1 in \cite{Bl01} 
that the multiplication in $A$ is separately weak$^*$-continuous, 
hence 
$A$ is actually a dual operator $G$-space.
Now we see that in either of the cases (i)~and~(ii) it follows by 
Theorem~3.1 that there exists a completely bounded idempotent 
mapping 
$P\colon A\to A$ with $\Ran P=A^G=\rho(G)'$, 
and we are done. 
\qed
\enddemo

We point out that some further smoothness properties of similarity orbits 
of group representations 
(in particular existence of complex structures on the unitary orbits) 
are discussed in \cite{Ma90}~and~\cite{MS95}.

\head References\endhead

\refstyle{ABCDEF}


\Refs\nofrills{}
\widestnumber\key{ABCDEF}

\ref\key{ACS95}
\by E.~Andruchow, G.~Corach, D.~Stojanoff 
\paper A geometric characterization of nuclearity and injectivity  
\jour J. Funct. Anal.  
vol 133  
\yr 1995
\pages no. 2, 474--494
\endref

\ref\key{AGG02}
\by A.~Arias, A.~Gheondea, S.~Gudder
\paper Fixed points of quantum operations
\jour J. Math. Phys.
\vol 43
\yr 2002
\pages no. 12, 5872--5881
\endref

\ref\key{BR04}
\by D.~Belti\c t\u a, T.S.~Ratiu 
\paper Symplectic leaves in real Banach Lie-Poisson spaces 
\jour Geom. Funct. Anal.
%\vol 
%\yr 
\pages (to appear). (See http://xxx.lanl.gov/abs/math.SG/0403345)
\endref

\ref\key{BH67}
\by J.F.~Berglund, K.H.~Hofmann
\book Compact Semitopological Semigroups and Weakly Almost Periodic Functions
\publ Lecture Notes in Mathematics, No. 42, Springer-Verlag
\publaddr Berlin-New York 
\yr 1967
\endref

\ref\key{BJM78}
\by J.F.~Berglund, H.D.~Junghenn, P.~Milnes
\book Compact Right Topological Semigroups and Generalizations of 
Almost Periodicity 
\publ Lecture Notes in Mathematics, No.~663, Springer-Verlag
\publaddr Berlin
\yr 1978
\endref

\ref\key{Bl01}
\by  D.P.~Blecher
\paper Multipliers and dual operator algebras
\jour J. Funct. Anal.
\vol 183
\yr 2001
\pages no. 2, 498--525
\endref


\ref\key{BL04}
\by D.P.~Blecher,  C.~Le Merdy
\book Operator Algebras and Their Modules. An Operator Space Approach
\publ London Mathematical Society Monographs. New Series.
The Clarendon Press, Oxford University Press
\publaddr Oxford
\yr 2004
\endref 

\ref\key{Bo65}
\by N.~Bourbaki 
\book \'El\'ements de Math\'ematique. Fasc. XIII. Livre VI: 
Int\'egration,  
Chapitres 1--4 %: In\'egalit\'es de convexit\'e, Espaces de Riesz, 
%Mesures sur les espaces localement compacts, Prolongement d'une mesure, 
%Espaces $L^p$ 
\publ (Deuxi\`eme \'edition revue et augment\'ee) 
Actualit\'es Scientifiques et Industrielles, No. 1175 Hermann
\publaddr Paris 
\yr 1965
\endref 

\ref\key{CG99}
\by G.~Corach, J.E.~Gal\'e 
\paper On amenability and geometry of spaces of bounded representations  
\jour J. London Math. Soc. (2)  
\vol 59  
\yr 1999
\pages no. 1, 311--329
\endref

\ref\key{DLRZ02}
\by S.~Doplicher, R.~Longo, J.E.~Roberts, L.~Zsid\'o 
\paper A remark on quantum group actions and nuclearity
\jour Rev. Math. Phys.  
\vol 14  
\yr 2002  
\pages no. 7-8, 787--796
\endref

\ref\key{ER00}
\by E.G.~Effros, Zh.-J.~Ruan
\book Operator Spaces
\publ London Mathematical Society Monographs. New Series, 23.
The Clarendon Press, Oxford University Press
\publaddr New York
\yr 2000
\endref

\ref\key{Gr69}
\by F.P.~Greenleaf
\book Invariant Means on Topological Groups and Their Applications
\publ Van Nostrand Mathematical Studies, No. 16, 
Van Nostrand Reinhold Co. 
\publaddr New York-Toronto, Ont.-London 
\yr 1969
\endref

\ref\key{dlH79}
\by P.~de la Harpe 
\paper Moyennabilit\'e du groupe unitaire et propri\'et\'e $P$ de 
Schwartz des alg\`ebres de von Neumann
\inbook in: Alg\`ebres d'Op\'erateurs (S\'em., Les Plans-sur-Bex, 1978)
\pages 220--227 
\publ Lecture Notes in Math., 725, Springer
\publaddr Berlin 
\yr 1979
\endref

\ref\key{HK77}
\by L.A.~Harris, W.~Kaup
\paper Linear algebraic groups in infinite dimensions
\jour Illinois J. Math. 
\vol 21 
\yr 1977 
\pages no. 3, 666--674
\endref

\ref\key{Ke02}
\by L.~K\'erchy
\paper Generalized Toeplitz operators
\jour Acta Sci. Math. (Szeged)
\vol 68
\yr 2002
\pages no. 1-2, 373--400
\endref

\ref\key{Ke04}
\by L.~K\'erchy
\paper Elementary and reflexive hyperplanes of generalized Toeplitz operators  
\jour J. Operator Theory  
\vol 51  
\yr 2004
\pages no. 2, 387--409
\endref 

\ref\key{Lu92}
\by A.~\L uczak
\paper Invariant states and ergodic dynamical systems on $W^*$-algebras
\jour Math. Proc. Cambridge Philos. Soc.  
\vol 111  
\yr 1992  
\pages no. 1, 181--192
\endref

\ref\key{Ma90}
\by M.~Martin 
\paper Projective representations of compact groups in $C^*$-algebras 
\inbook in: Linear Operators in Function Spaces (Timi\c soara, 1988) 
\pages 237--253 
\publ Oper. Theory Adv. Appl., Birkh\"auser 
\vol 43
\publaddr Basel
\yr 1990
\endref

\ref\key{MS95}
\by M.~Martin, N.~Salinas 
\paper Differential geometry of generalized Grassmann manifolds in $C^*$-algebras  
\inbook in: Operator Theory and Boundary Eigenvalue Problems (Vienna, 1993)
\pages 206--243 
\publ Oper. Theory Adv. Appl., Birkh\"auser
\vol 80
\publaddr Basel
\yr 1995
\endref

\ref\key{OR03}
\by A.~Odzijewicz, T.S.~Ratiu 
\paper Banach Lie-Poisson spaces and reduction
\jour Comm. Math. Phys.  
\vol 243  
\yr 2003
\pages   no.~1, 1--54
\endref

\ref\key{Pat92}
\by A.L.T.~Paterson 
\paper Nuclear $C^*$-algebras have amenable unitary groups  
\jour Proc. Amer. Math. Soc.  
\vol 114  
\yr 1992
\pages no. 3, 719--721
\endref 

\ref\key{Pa02}
\by V.~Paulsen 
\book Completely Bounded Maps and Operator Algebras
\publ Cambridge Studies in Advanced Mathematics, 78. Cambridge University Press 
\publaddr Cambridge
\yr 2002
\endref 

\ref\key{Pe79}
\by G.K.~Pedersen 
\book $C^*$-algebras and Their Automorphism Groups  
\publ London Mathematical Society Monographs, 14. Academic Press, Inc. 
\publaddr London-New York
\yr 1979
\endref

\ref\key{Po03}
\by G.~Popescu
\paper Similarity and ergodic theory of positive linear maps
\jour J. Reine Angew. Math. 
\vol 561 
\yr 2003
\pages 87--129
\endref

\ref\key{Sch66}
\by H.H.~ Schaefer 
\book Topological Vector Spaces 
\publ The Macmillan Co., Collier-Macmillan Ltd.
\publaddr New York-London
\yr 1966
\endref

\ref\key{SS98}
\by A.M.~Sinclair, R.R.~Smith
\paper The Hochschild cohomology problem for von Neumann algebras  
\jour Proc. Natl. Acad. Sci. USA  
\vol 95  
\yr 1998
\pages no. 7, 3376--3379
\endref

\ref\key{SS04}
\by A.M.~Sinclair, R.R.~Smith
\paper A survey of Hochschild cohomology for von Neumann algebras
\jour Contemp. Math 
%\vol 
%\yr 
\pages (to appear)%. (See http://www.math.tamu.edu/~roger.smith/preprints.html.)
\endref

\ref\key{Ta03}
\by M.~Takesaki
\book Theory of Operator Algebras III
\publ Encyclopaedia of Mathematical Sciences, 127. 
Operator Algebras and Non-commutative Geometry, 8. Springer-Verlag 
\publaddr Berlin 
\yr 2003
\endref

\ref\key{Up85}
\by H.~Upmeier
\book Symmetric Banach Manifolds and Jordan $C^*$-algebras 
\publ North-Holland Mathematics Studies, 104. 
Notas de Matem\'atica, 96. North-Holland Publishing Co. 
\publaddr Amsterdam 
\yr 1985
\endref

\ref\key{Wa72}
\by G.~Warner 
\book Harmonic Analysis on Semi-Simple Lie Groups I 
\publ Die Grundlehren der mathematischen Wissenschaften, Band 188. 
Springer-Verlag 
\publaddr New York-Heidelberg
\yr 1972
\endref



\endRefs
\enddocument